\documentclass[a4paper]{article}
\usepackage{amssymb}
\usepackage{graphicx}
\usepackage{hyperref}
\usepackage{latexsym}
\newtheorem{theorem}{Theorem}
\newtheorem{corollary}{Corollary}
\newtheorem{lemma}{Lemma}
\newtheorem{definition}{Definition}

\newtheorem{query}{Query}

\begin{document}
\title{A case for weakening the Church-Turing Thesis}
\author{\normalsize{Bhupinder Singh Anand}}
\date{\small{Draft of \today.\footnote{Subject class: LO; MSC: 03B10}}}
\maketitle

\begin{abstract}We conclude from G\"{o}del's Theorem VII of his seminal 1931 paper that every recursive function $f(x_{1}, x_{2})$ is representable in the first-order Peano Arithmetic PA by a formula $[F(x_{1}, x_{2}, x_{3})]$ which is algorithmically verifiable, but not algorithmically computable, \textit{if} we assume that the negation of a universally quantified formula of the first-order predicate calculus is always indicative of the existence of a counter-example under the standard interpretation of PA. We conclude that the standard postulation of the Church-Turing Thesis does not hold if we define a number-theoretic formula as effectively computable if, and only if, it is algorithmically verifiable; and needs to be replaced by a weaker postulation of the Thesis as an equivalence.

\vspace{+1ex}
\end{abstract}

\section{Introduction}
\label{sec:2.1}

We begin by noting that the following theses are classically equivalent\footnote{cf.\ \cite{Me64}, p.237.}:

\vspace{+1ex}
\textbf{Standard Church's Thesis}\footnote{\textit{Church's (original) Thesis} The effectively computable number-theoretic functions are the algorithmically computable number-theoretic functions \cite{Ch36}.} A number-theoretic function (or relation, treated as a Boolean function) is effectively computable if, and only if, it is partial-recursive\footnote{cf.\ \cite{Me64}, p.227.}.

\vspace{+1ex}
\textbf{Standard Turing's Thesis}\footnote{After describing what he meant by ``computable" numbers in the opening sentence of his 1936 paper on Computable Numbers \cite{Tu36}, Turing immediately expressed this thesis---albeit informally---as: ``\ldots the computable numbers include all numbers which could naturally be regarded as computable".} A number-theoretic function (or relation, treated as a Boolean function) is effectively computable if, and only if, it is Turing-computable\footnote{cf.\ \cite{BBJ03}, p.33.}.

\vspace{+1ex}
In this paper we shall argue that the principle of Occam's razor suggests the Theses should be postulated minimally as the following equivalences:

\vspace{+1ex}
\textbf{Weak Church's Thesis} A number-theoretic function (or relation, treated as a Boolean function) is effectively computable if, and only if, it is instantiationally equivalent to a partial-recursive function (or relation, treated as a Boolean function).

\vspace{+1ex}
\textbf{Weak Turing's Thesis} A number-theoretic function (or relation, treated as a Boolean function) is effectively computable if, and only if, it is instantiationally equivalent to a Turing-computable function (or relation, treated as a Boolean function).

\subsection{The need for explicitly distinguishing between `instantiational' and `uniform' methods}
\label{sec:9}

It is significant that both Kurt G\"{o}del (initially) and Alonzo Church (subseque- ntly---possibly under the influence of G\"{o}del's disquietitude) enunciated Church's formulation of `effective computability' as a Thesis because G\"{o}del was instinctively uncomfortable with accepting it as a definition that \textit{minimally} captures the essence of `\textit{intuitive} effective computability'\footnote{See \cite{Si97}.}.

G\"{o}del's reservations seem vindicated if we accept that a number-theoretic function can be effectively computable instantiationally (in the sense of being algorithmically \textit{verifiable} as defined in Section \ref{sec:1.02} below), but not by a uniform method (in the sense of being algorithmically \textit{computable} as defined in Section \ref{sec:1.02}).

The significance of the fact (considered below in Section \ref{sec:5.4.0}) that `truth' too can be effectively decidable \textit{both} instantiationally \textit{and} by a uniform (algorithmic) method under the standard interpretation of PA is reflected in G\"{o}del's famous 1951 Gibbs lecture\footnote{\cite{Go51}.}, where he remarks: 

\begin{quote}
``I wish to point out that one may conjecture the truth of a universal proposition (for example, that I shall be able to verify a certain property for any integer given to me) and at the same time conjecture that no general proof for this fact exists. It is easy to imagine situations in which both these conjectures would be very well founded. For the first half of it, this would, for example, be the case if the proposition in question were some equation $F(n) = G(n)$ of two number-theoretical functions which could be verified up to very great numbers $n$."\footnote{Parikh's paper \cite{Pa71} can also be viewed as an attempt to investigate the consequences of expressing the essence of G\"{o}del's remarks formally.} 
\end{quote}

Such a possibility is also implicit in Turing's remarks\footnote{\cite{Tu36}, \S9(II), p.139.}:

\begin{quote}
``The computable numbers do not include all (in the ordinary sense) definable numbers. Let P be a sequence whose \textit{n}-th figure is 1 or 0 according as \textit{n} is or is not satisfactory. It is an immediate consequence of the theorem of \S8 that P is not computable. It is (so far as we know at present) possible that any assigned number of figures of P can be calculated, but not by a uniform process. When sufficiently many figures of P have been calculated, an essentially new method is necessary in order to obtain more figures." 
\end{quote}

The need for placing such a distinction on a formal basis has also been expressed explicitly on occasion\footnote{Parikh's distinction between `decidability' and `feasibility' in \cite{Pa71} also appears to echo the need for such a distinction.}. Thus, Boolos, Burgess and Jeffrey\footnote{\cite{BBJ03}, p. 37.} define a diagonal function, $d$, any value of which can be decided effectively, although there is no single algorithm that can effectively compute $d$. 

Now, the straightforward way of expressing this phenomenon should be to say that there are well-defined number-theoretic functions that are effectively computable instantiationally but not algorithmically. Yet, following Church and Turing, such functions are labeled as uncomputable\footnote{The issue here seems to be that, when using language to express the abstract objects of our individual, and common, mental `concept spaces', we use the word `exists' loosely in three senses, without making explicit distinctions between them (see \cite{An07c}).}!

\begin{quote}
\noindent ``According to Turing's Thesis, since $d$ is not Turing-computable, $d$ cannot be effectively computable. Why not? After all, although no Turing machine computes the function $d$, we were able to compute at least its first few values, For since, as we have noted, $f_{1} = f_{1} = f_{1} =$ the empty function we have $d(1) = d(2) = d(3) = 1$. And it may seem that we can actually compute $d(n)$ for any positive integer $n$---if we don't run out of time."\footnote{\cite{BBJ03}, p.37.}
\end{quote}

The reluctance to treat a function such as $d(n)$---or the function $\Omega(n)$ that computes the $n^{th}$ digit in the decimal expression of a Chaitin constant $\Omega$\footnote{Chaitin's Halting Probability is given by $0 < \Omega = \sum2^{-|p|} < 1$, where the summation is over all self-delimiting programs $p$ that halt, and $|p|$ is the size in bits of the halting program $p$; see \cite{Ct75}.}---as computable, on the grounds that the `time' needed to compute it increases monotonically with $n$, is curious\footnote{The incongruity of this is addressed by Parikh in \cite{Pa71}.}; the same applies to any total Turing-computable function $f(n)$\footnote{The only difference being that, in the latter case, we know there is a common `program' of constant length that will compute $f(n)$ for any given natural number $n$; in the former, we know we may need distinctly different programs for computing $f(n)$ for different values of $n$, where the length of the program will, sometime, reference $n$.}!

\subsection{Distinguishing between algorithmically verifiability and algorithmic computability}
\label{sec:2.1.0}

We now show in Theorem \ref{sec:2.3.thm.1} that if Aristotle's particularisation\footnote{See Definition \ref{def:A1} below.} is presumed valid over the structure $\mathbb{N}$ of the natural numbers---as is the case under the standard interpretation of the first-order Peano Arithmetic PA---then it follows from the instantiational nature of the (constructively defined\footnote{By Kurt G\"{o}del; see Appendix A, Section \ref{sec:5.4.1}, Lemma 1}) G\"{o}del $\beta$-function that a primitive recursive relation can be instantiationally equivalent to an arithmetical relation, where the former is algorithmically computable over $\mathbb{N}$, whilst the latter is algorithmically verifiable but not algorithmically computable over $\mathbb{N}$.

\begin{quote}
\footnotesize
\noindent \textbf{Analagous distinctions in analysis:} The distinction between algorithmically computable, and algorithmically verifiable but not algorithmically computable, number-theoretic functions seeks to reflect in arithmetic the essence of \textit{uniform} methods\footnote{See Section \ref{sec:9}.}, classically characterised by the distinctions in analysis between: (a) uniformly continuous, and point-wise continuous but not uniformly continuous, functions over an interval; (b) uniformly convergent, and point-wise convergent but not uniformly convergent, series.

\vspace{+1ex}
\noindent \textbf{A limitation of set theory and a possible barrier to computation:} We note, further, that the above distinction cannot be reflected within a language---such as the set theory ZF---which identifies `equality' with `equivalence'. Since functions are defined extensionally as mappings, such a language cannot recognise that a set which represents a primitive recursive function may be equivalent to, but computationally different from, a set that represents an arithmetical function; where the former function is algorithmically computable over $\mathbb{N}$, whilst the latter is algorithmically verifiable but not algorithmically computable over $\mathbb{N}$.
\end{quote}

\subsubsection{Significance of G\"{o}del's $\beta$-function}
\label{sec:2.1.0}

In Theorem VII \footnote{\cite{Go31}, pp.30-31; reproduced below in Appendix A, Section \ref{sec:5.4.1}.} of his seminal 1931 paper on formally undecidable arithmetical propositions, G\"{o}del showed that, given a total number-theoretic function $f(x)$ and any natural number $n$, we can construct a primitive recursive function $\beta(z, y, x)$ and natural numbers $b_{n}, c_{n}$ such that $\beta(b_{n}, c_{n}, i)$ $= f(i)$ for all $0 \leq i \leq n$.

In this paper we shall essentially answer the following question affirmatively:

\begin{query}
\noindent Does G\"{o}del's Theorem VII admit construction of an arithmetical function $A(x)$ such that:

\begin{quote}
\noindent (a) for any given natural number $n$, there is an algorithm that can verify $A(i) = f(i)$ for all $0 \leq i \leq n$ (hence $A(x)$ may be said to be algorithmically verifiable if $f(x)$ is recursive);

\vspace{+1ex}
\noindent (b) there is no algorithm that can verify $A(i) = f(i)$ for all $0 \leq i$ (so $A(x)$ may be said to be algorithmically uncomputable)?
\end{quote}
\end{query}

\subsubsection{Defining effective computability}
\label{sec:2.1.0.1}

We shall then formally define what it means for a formula of an arithmetical language to be:

\begin{quote}
(i) Algorithmically verifiable;

\vspace{+.5ex}
(ii) Algorithmically computable.
\end{quote}

\noindent under an interpretation.

\vspace{+1ex}
We shall further propose the definition:

\vspace{+1ex}
\noindent \textbf{Effective computability:} A number-theoretic formula is effectively computable if, and only if, it is algorithmically verifiable. 

\begin{quote}
\footnotesize
\textbf{Intuitionistically unobjectionable:} We note first that since every finite set of integers is recursive, every well-defined number-theoretical formula is algorithmically verifiable, and so the above definition is intuitionistically unobjectionable; and second that the existence of an arithmetic formula that is algorithmically verifiable but not algorithmically computable (Theorem \ref{sec:2.3.thm.1}) supports G\"{o}del's reservations on Alonzo Church's original intention to label his Thesis as a definition\footnote{See Section \ref{sec:9}.}.
\end{quote}

We shall then show that the algorithmically verifiable and the algorithmically computable PA formulas are well-defined under the standard interpretation of PA since:

\begin{quote}
(a) The PA-formulas are decidable as satisfied / unsatisfied or true / false under the standard interpretation of PA if, and only if, they are algorithmically verifiable;

\vspace{+.5ex}
(b) The algorithmically computable PA-formulas are a proper subset of the algorithmically verifiable PA-formulas;

\vspace{+,5ex}
(c) The PA-axioms are algorithmically computable as satisfied / true under the standard interpretation of PA;

\vspace{+.5ex}
(d) Generalisation and Modus Ponens preserve algorithmically computable truth under the standard interpretation of PA;

\vspace{+.5ex}
(e) The provable PA-formulas are precisely the ones that are algorithmically computable as satisfied / true under the standard interpretation of PA.
\end{quote}

\section{Comments, Notation and Standard Definitions}
\label{sec:A}

\noindent \textbf{Comments} We have taken some liberty in emphasising standard definitions selectively, and interspersing our arguments liberally with comments and references, generally of a foundational nature. These are intended to reflect our underlying thesis that essentially arithmetical problems appear more natural when expressed---and viewed---within the perspective of an interpretation of PA that appeals to the \textit{evidence} provided by a deterministic Turing machine\footnote{A deterministic Turing machine has only one possible move from a given configuration.} along the lines suggested in Section \ref{sec:5.4.0}; a perspective that, by its very nature, cannot appeal implicitly to transfinite concepts.

\begin{quote}
\footnotesize \textbf{Evidence} ``It is by now folklore \ldots that one can view the \textit{values} of a simple functional language as specifying \textit{evidence} for propositions in a constructive logic \ldots"\footnote{\cite{Mu91}.}.
\end{quote}

\vspace{+1ex}
\noindent \textbf{Notation} We use square brackets to indicate that the contents represent a symbol or a formula---of a formal theory---generally assumed to be well-formed unless otherwise indicated by the context.\footnote{In other words, expressions inside the square brackets are to be only viewed syntactically as juxtaposition of symbols that are to be formed and manipulated upon strictly in accordance with specific rules for such formation and manipulation---in the manner of a mechanical or electronic device---without any regards to what the symbolism might represent semantically under an interpretation that gives them meaning.} We use an asterisk to indicate that the expression is to be interpreted semantically with respect to some well-defined interpretation.

\begin{definition}
\label{def:A1}
\noindent \textbf{Aristotle's particularisation} This holds that from an assertion such as:

\begin{quote}
`It is not the case that: for any given $x$, $P^{*}(x)$ does not hold',
\end{quote}

\noindent usually denoted symbolically by `$\neg(\forall x)\neg P^{*}(x)$', we may always validly infer in the classical, Aristotlean, logic of predicates\footnote{\cite{HA28}, pp.58-59.} that:

\begin{quote}
`There exists an unspecified $x$ such that $P^{*}(x)$ holds',
\end{quote}

\noindent usually denoted symbolically by `$(\exists x)P^{*}(x)$'.
\end{definition}

\begin{quote}
\footnotesize
\textbf{The significance of Aristotle's particularisation for the first-order predicate calculus:} We note that in a formal language the formula `$[(\exists x)P(x)]$' is an abbreviation for the formula `$[\neg(\forall x)\neg P(x)]$'. The commonly accepted interpretation of this formula---and a fundamental tenet of classical logic unrestrictedly adopted as intuitively obvious by standard literature that seeks to build upon the formal first-order predicate calculus\footnote{See \cite{Hi25}, p.382; \cite{HA28}, p.48; \cite{Sk28}, p.515; \cite{Go31}, p.32.; \cite{Kl52}, p.169; \cite{Ro53}, p.90; \cite{BF58}, p.46; \cite{Be59}, pp.178 \& 218; \cite{Su60}, p.3; \cite{Wa63}, p.314-315; \cite{Qu63}, pp.12-13; \cite{Kn63}, p.60; \cite{Co66}, p.4; \cite{Me64}, p.52(ii); \cite{Nv64}, p.92; \cite{Li64}, p.33; \cite{Sh67}, p.13; \cite{Da82}, p.xxv; \cite{Rg87}, p.xvii; \cite{EC89}, p.174; \cite{Mu91}; \cite{Sm92}, p.18, Ex.3; \cite{BBJ03}, p.102.}---tacitly appeals to Aristotlean particularisation.

However, L. E. J. Brouwer had noted in his seminal 1908 paper on the unreliability of logical principles\footnote{\cite{Br08}.} that the commonly accepted interpretation of this formula is ambiguous if interpretation is intended over an infinite domain. He essentially argued that, even supposing the formula `$[P(x)]$' of a formal Arithmetical language interprets as an arithmetical relation denoted by `$P^{*}(x)$', and the formula `$[\neg(\forall x)\neg P(x)]$' as the arithmetical proposition denoted by `$\neg(\forall x)\neg P^{*}(x)$', the formula `$[(\exists x)P(x)]$' need not interpret as the arithmetical proposition denoted by the usual abbreviation `$(\exists x)P^{*}(x)$'; and that such postulation is invalid as a general logical principle in the absence of a means for constructing some putative object $a$ for which the proposition $P^{*}(a)$ holds in the domain of the interpretation.

Hence we shall follow the convention that the assumption that `$(\exists x)P^{*}(x)$' is the intended interpretation of the formula `$[(\exists x)P(x)]$'---which is essentially the assumption that Aristotle's particularisation holds over the domain of the interpretation---must always be explicit.

\vspace{+1ex}
\noindent \textbf{The significance of Aristotle's particularisation for PA:} In order to avoid intuitionistic objections to his reasoning, G\"{o}del introduced the syntactic property of $\omega$-consistency as an explicit assumption in his formal reasoning in his seminal 1931 paper on formally undecidable arithmetical propositions\footnote{\cite{Go31}, p.23 and p.28.}.

G\"{o}del explained at some length\footnote{In his introduction on p.9 of \cite{Go31}.} that his reasons for introducing $\omega$-consistency explicitly was to avoid appealing to the semantic concept of classical arithmetical truth in Aristotle's logic of predicates (which presumes Aristotle's particularisation).

However, we note that the two concepts are meta-mathematically equivalent in the sense that, if PA is consistent, then PA is $\omega$-consistent if, and only if, Aristotle's particularisation holds under the standard interpretation of PA.
\end{quote}

\begin{definition}
\label{def:A2}
\noindent \textbf{The structure $\mathbb{N}$} The structure of the natural numbers---namely, \{$N$ (\textit{the set of natural numbers}); $=$ (\textit{equality}); $'$ (\textit{the successor function}); $+$ (\textit{the addition function}); $ \ast $ (\textit{the product function}); $0$ (\textit{the null} \textit{element})\}.
\end{definition}

\begin{definition}
\label{def:A3}
\noindent \textbf{The axioms of first-order Peano Arithmetic (PA)}
\end{definition}

\begin{tabbing}
\textbf{PA$_{1}$} \= $[(x_{1} = x_{2}) \rightarrow ((x_{1} = x_{3}) \rightarrow (x_{2} = x_{3}))]$; \\

\textbf{PA$_{2}$} \> $[(x_{1} = x_{2}) \rightarrow (x_{1}^{\prime} = x_{2}^{\prime})]$; \\

\textbf{PA$_{3}$} \> $[0 \neq x_{1}^{\prime}]$; \\

\textbf{PA$_{4}$} \> $[(x_{1}^{\prime} = x_{2}^{\prime}) \rightarrow (x_{1} = x_{2})]$; \\

\textbf{PA$_{5}$} \> $[( x_{1} + 0) = x_{1}]$; \\

\textbf{PA$_{6}$} \> $[(x_{1} + x_{2}^{\prime}) = (x_{1} + x_{2})^{\prime}]$; \\

\textbf{PA$_{7}$} \> $[( x_{1} \star 0) = 0]$; \\

\textbf{PA$_{8}$} \> $[( x_{1} \star x_{2}^{\prime}) = ((x_{1} \star x_{2}) + x_{1})]$; \\

\textbf{PA$_{9}$} \> For any well-formed formula $[F(x)]$ of PA: \\

\> $[F(0) \rightarrow (((\forall x)(F(x) \rightarrow F(x^{\prime}))) \rightarrow (\forall x)F(x))]$.
\end{tabbing}

\begin{definition}
\label{def:A4}
\noindent \textbf{Generalisation in PA} If $[A]$ is PA-provable, then so is $[(\forall x)A]$.
\end{definition}

\begin{definition}
\label{def:A5}
\noindent \textbf{Modus Ponens in PA} If $[A]$ and $[A \rightarrow B]$ are PA-provable, then so is $[B]$.
\end{definition}

\begin{definition}
\label{def:A6}
\noindent \textbf{Standard interpretation of PA} The standard interpretation of PA over the structure $\mathbb{N}$ is the one in which the logical constants have their `usual' interpretations\footnote{See \cite{Me64}, p.49.} in Aristotle's logic of predicates (which subsumes Aristotle's particularisation), and\footnote{See \cite{Me64}, p.107.}:

\begin{tabbing}
(a) \= the set of non-negative integers is the domain; \\

(b) \> the symbol [0] interprets as the integer 0; \\

(c) \> the symbol $[']$ interprets as the successor operation (addition of 1); \\

(d) \> the symbols $[+]$ and $[*]$ interpret as ordinary addition and multiplication; \\

(e) \> the symbol $[=]$ interprets as the identity relation.
\end{tabbing}
\end{definition}

\section{G\"{o}del's Theorem VII and algorithmically verifiable, but not algorithmically computable, arithmetical propositions}
\label{sec:2.1.1}

In his seminal 1931 paper on formally undecidable arithmetical propositions, G\"{o}del defined a curious primitive recursive function---G\"{o}del's $\beta$-function---as\footnote{cf.\ \cite{Go31}, p.31, Lemma 1; \cite{Me64}, p.131, Proposition 3.21.}:

\begin{definition}
\label{sec:2.2.1.def.0}
$\beta (x_{1}, x_{2}, x_{3}) = rm(1+(x_{3}+ 1) \star x_{2}, x_{1})$

\vspace{+.5ex}
\noindent where $rm(x_{1}, x_{2})$ denotes the remainder obtained on dividing $x_{2}$ by $x_{1}$.
\end{definition}

\vspace{+1ex}
G\"{o}del showed that the above function has the remarkable property that:

\begin{lemma}
\label{sec:2.2.2.lem.1}
For any given denumerable sequence of natural numbers, say $f(k, 0), f(k, 1), \ldots$, and any given natural number $n$, we can construct natural numbers $b, c, j$ such that:

\begin{quote}
(i) $j = max(n, f(k, 0), f(k, 1), \ldots, f(k, n))$;

(ii) $c = j$!;

(iii) $\beta(b, c, i) = f(k, i)$ for $0 \leq i \leq n$.
\end{quote}
\end{lemma}

\vspace{+1ex}
\noindent \textbf{Proof} This is a standard result\footnote{cf.\ \cite{Go31}, p.31, p.31, Lemma 1; \cite{Me64}, p.131, Proposition 3.22.}.\hfill $\Box$

\begin{quote}
\footnotesize
\textbf{G\"{o}del's original argument:} We reproduce G\"{o}del's original argument---which yields this critical lemma---in Appendix A, Section \ref{sec:5.4.1}.
\end{quote}

Now we have the standard definition\footnote{\cite{Me64}, p.118.}:

\begin{definition}
\label{sec:2.2.1.def.1}
A number-theoretic function $f(x_{1}, \ldots, x_{n})$ is said to be representable in PA if, and only if, there is a PA formula $[F(x_{1}, \dots, x_{n+1})]$ with the free variables $[x_{1}, \ldots, x_{n+1}]$, such that, for any given natural numbers $k_{1}, \ldots, k_{n+1}$:

\begin{quote}
(i) if $f(k_{1}, \ldots, k_{n}) = k_{n+1}$ then PA proves: $[F(k_{1}, \ldots, k_{n}, k_{n+1})]$;

(ii) PA proves: $[(\exists_{1} x_{n+1})F(k_{1}, \ldots, k_{n}, x_{n+1})]$.
\end{quote}

The function $f(x_{1}, \ldots, x_{n})$ is said to be strongly representable in PA if we further have that:

\begin{quote}
(iii) PA proves: $[(\exists_{1} x_{n+1})F(x_{1}, \ldots, x_{n}, x_{n+1})]$
\end{quote}
\end{definition}

\begin{quote}
\footnotesize
\textbf{Interpretation of `$[\exists_{1}]$':} The symbol `$[\exists_{1}]$' denotes `uniqueness' under an interpretation which assumes that Aristotle's particularisation holds in the domain of the interpretation. Formally, however, the PA formula $[(\exists_{1} x_{3})F(x_{1}, x_{2}, x_{3})]$ is merely a short-hand notation for the PA formula $[\neg(\forall x_{3})\neg F(x_{1}, x_{2}, x_{3}) \wedge (\forall y)(\forall z)$ $(F(x_{1}, x_{2}, y) \wedge F(x_{1}, x_{2}, z) \rightarrow y=z)]$.
\end{quote}

We then have:

\begin{lemma}
\label{sec:2.2.1.lem.1}
$\beta(x_{1}, x_{2}, x_{3})$ is strongly represented in PA by $[Bt(x_{1}, x_{2}, x_{3}, x_{4})]$, which is defined as follows: 

\begin{quote}
$[(\exists w)(x_{1} = ((1 + (x_{3} + 1)\star x_{2}) \star w + x_{4}) \wedge (x_{4} < 1 + (x_{3} + 1) \star x_{2}))]$.
\end{quote}
\end{lemma}

\vspace{+1ex}
\noindent \textbf{Proof} This is a standard result\footnote{cf. \cite{Me64}, p.131, proposition 3.21.}.\hfill $\Box$

\vspace{+1ex}
G\"{o}del further showed (also under the tacit, but critical, presumption of Aristotle's particularisation\footnote{The implicit assumption being that the negation of a universally quantified formula of the first-order predicate calculus is indicative of ``the existence of a counter-example"---\cite{Go31}, p.32.}) that:

\begin{lemma}
\label{sec:2.2.3}
If $f(x_{1}, x_{2})$ is a recursive function defined by:

\begin{quote}
(i) $f(x_{1}, 0) = g(x_{1})$

(ii) $f(x_{1}, (x_{2}+1)) = h(x_{1}, x_{2}, f(x_{1}, x_{2}))$
\end{quote}

\noindent where $g(x_{1})$ and $h(x_{1}, x_{2}, x_{3})$ are recursive functions of lower rank\footnote{cf.\ \cite{Me64}, p.132; \cite{Go31}, p.30(2).} that are represented in PA by well-formed formulas $[G(x_{1}, x_{2})]$ and $[H(x_{1}, x_{2}, x_{3}, x_{4})]$, then $f(x_{1}, x_{2})$ is represented in PA by the following well-formed formula, denoted by $[F(x_{1}, x_{2}, x_{3})]$:

\begin{quote}
$[(\exists u)(\exists v)(((\exists w)(Bt(u, v, 0, w) \wedge G(x_{1}, w))) \wedge Bt(u, v, x_{2}, x_{3}) \wedge (\forall w)(w$ $< x_{2} \rightarrow (\exists y)(\exists z)(Bt(u, v, w, y) \wedge Bt(u, v, (w+1), z) \wedge H(x_{1}, w, y, z)))]$.
\end{quote}
\end{lemma}

\vspace{+1ex}
\noindent \textbf{Proof} This is a standard result\footnote{cf.\ \cite{Go31}, p.31(2); \cite{Me64}, p.132.}.\hfill $\Box$

\begin{quote}
\footnotesize
\textbf{G\"{o}del's original argument:} We reproduce G\"{o}del's original argument and proof of this critical lemma in Appendix A, Section \ref{sec:5.4.1}.
\end{quote}

\subsection{What does ``$[(\exists_{1} x_{3})F(k, m, x_{3})]$ is provable" assert under the standard interpretation of PA?}

Now, if the PA formula $[F(x_{1}, x_{2}, x_{3})]$ represents in PA the recursive function denoted by $f(x_{1}, x_{2})$ then by definition, for any given numerals $[k], [m]$, the formula $[(\exists_{1} x_{3})F(k, m, x_{3})]$ is provable in PA; and true under the standard interpretation of PA. We thus have that:

\begin{lemma}
\label{sec:2.2.3.lem.1}

``$[(\exists_{1} x_{3})F(k,$ $m, x_{3})]$ is true under the standard interpretation of PA" is the assertion that:

\begin{quote}
Given any natural numbers $k, m$, we can construct natural numbers $t_{(k, m)},$ $u_{(k, m)},$ $v_{(k, m)}$---all functions of $k, m$---such that:

\vspace{+1ex}
\begin{quote}

(a) $\beta(u_{(k, m)}, v_{(k, m)}, 0) = g(k)$;

\vspace{+.5ex}
(b) for all $i<m$, $\beta(u_{(k, m)}, v_{(k, m)}, i) = h(k, i, f(k, i))$;

\vspace{+.5ex}
(c) $\beta(u_{(k, m)}, v_{(k, m)}, m) = t_{(k, m)}$;
 
\end{quote}

\vspace{+1ex}
where $f(x_{1}, x_{2})$, $g(x_{1})$ and $h(x_{1}, x_{2}, x_{3})$ are any recursive functions that are formally represented in PA by $F(x_{1}, x_{2},$ $x_{3}), G(x_{1}, x_{2})$ and $H(x_{1}, x_{2}, x_{3},$ $x_{4})$ respectively such that:

\vspace{+1ex}
\begin{quote}
(i)      $f(k, 0) = g(k)$

\vspace{+.5ex}
(ii)     $f(k, (y+1)) = h(k, y, f(k, y))$ for all $y<m$

\vspace{+.5ex}
(iii)    $g(x_{1})$ and $h(x_{1}, x_{2}, x_{3})$ are recursive functions that are assumed to be of lower rank than $f(x_{1}, x_{2})$.
\end{quote}
\end{quote}
\end{lemma}

\noindent \textbf{Proof} For any given natural numbers $k$ and $m$, if $[F(x_{1}, x_{2}, x_{3})]$ interprets as a well-defined arithmetical relation under the standard interpretation of PA, then we can define a deterministic Turing machine $TM$ that can `construct' the sequences $f(k, 0), f(k, 1), \ldots, f(k, m)$ and $\beta(u_{(k, m)}, v_{(k, m)}, 0),$ $\beta(u_{(k, m)}, v_{(k, m)}, 1),$ $\ldots, \beta(u_{(k, m)}, v_{(k, m)}, m)$ and give evidence to verify the assertion.\hfill $\Box$

\begin{quote}
\footnotesize
\textbf{Does $[F(x_{1}, x_{2}, x_{3}]$ interpret as a well-defined predicate?} A critical issue that we do not address here is whether the PA formula $[F(x_{1}, x_{2}, x_{3}]$ can be considered to interpret under a sound interpretation of PA as a well-defined predicate, since the denumerable sequences $\{f(k, 0), f(k, 1),$ $\ldots , f(k, m), m_{p}: p>0$ and $m_{p}$ is not equal to $m_{q}$ if $p$ is not equal to $q\}$---are represented by denumerable, distinctly different, functions $\beta(u_{p_{1}}, v_{p_{2}}, i)$ respectively. There are thus denumerable pairs $(u_{p_{1}}, v_{p_{2}})$ for which $\beta(u_{p_{1}}, v_{p_{2}}, i)$ yields the sequence $f(k, 0), f(k, 1),$ $\ldots , f(k, m)$.
\end{quote}

\vspace{1ex}
We now see that:

\begin{theorem}
\label{sec:2.3.thm.1}
Under the standard interpretation of PA $[(\exists_{1} x_{3})F(x_{1}, x_{2},$ $x_{3})]$ is algorithmically verifiable, but not algorithmically computable, as always true over $\mathbb{N}$.
\end{theorem}

\vspace{+1ex}
\noindent \textbf{Proof} It follows from Lemma \ref{sec:2.2.3.lem.1} that:

\begin{quote}
(1) $[(\exists_{1} x_{3})F(k, m, x_{3})]$ is PA-provable for any given numerals $[k, m]$. Hence $[(\exists_{1} x_{3})F(k, m, x_{3})]$ is true under the standard interpretation of PA. It then follows from the definition of $[F(x_{1}, x_{2}, x_{3})]$ in Lemma \ref{sec:2.2.3} that, for any given natural numbers $k, m$, we can construct some pair of natural numbers $u_{(k, m)}, v_{(k, m)}$---where $u_{(k, m)}, v_{(k, m)}$ are functions of the given natural numbers $k$ and $m$---such that:

\vspace{+1ex}
\begin{quote}
(a) $\beta(u_{(k, m)}, v_{(k, m)}, i) = f(k, i)$ for $0 \leq i \leq m$;

\vspace{+.5ex}
(b) $F^{*}(k, m, f(k, m))$ holds in $\mathbb{N}$.
\end{quote}

\vspace{+1ex}
Since $\beta(x_{1}, x_{2}, x_{3})$ is primitive recursive, $\beta(u_{(k, m)}, v_{(k, m)}, i)$ defines a deterministic Turing machine $TM$ that can `construct' the denumerable sequence $f'(k, 0), f'(k, 1),$ $\ldots$ for any given natural numbers $k$ and $m$ such that:

\vspace{+1ex}
\begin{quote}
(c) $f(k, i) = f'(k, i)$ for $0 \leq i \leq m$.
\end{quote}

\vspace{+1ex}
We can thus define a deterministic Turing machine $TM$ that will give evidence that the PA formula $[(\exists_{1} x_{3})F(k, m, x_{3})]$ is true under the standard interpretation of PA.

\vspace{+.5ex}
Hence $[(\exists_{1} x_{3})F(x_{1}, x_{2}, x_{3})]$ is algorithmically verifiable over $\mathbb{N}$ under the standard interpretation of PA.

\vspace{+1ex}
(2) Now, the pair of natural numbers $u_{(x_{1}, x_{2})}, v_{(x_{1}, x_{2})}$ are defined such that:

\vspace{+1ex}
\begin{quote}
(a) $\beta(u_{(x_{1}, x_{2})}, v_{(x_{1}, x_{2})}, i) = f(x_{1}, i)$ for $0 \leq i \leq x_{2}$;

\vspace{+.5ex}
(b) $F^{*}(x_{1}, x_{2}, f(x_{1}, x_{2}))$ holds in $\mathbb{N}$;
\end{quote}

\vspace{+1ex}
where $v_{(x_{1}, x_{2})}$ is defined in Lemma \ref{sec:2.2.3} as $j$!, and:

\vspace{+1ex}
\begin{quote}
(c) $j = max(n, f(x_{1}, 0), f(x_{1}, 1), \ldots, f(x_{1}, x_{2}))$;

\vspace{+.5ex}
(d) $n$ is the `number' of terms in the sequence $f(x_{1}, 0), f(x_{1}, 1),$ $\ldots, f(x_{1}, x_{2})$.
\end{quote}

\vspace{+1ex}
Since $j$ is not definable for a denumerable sequence $\beta(u_{(x_{1}, x_{2})},$ $v_{(x_{1}, x_{2})},$ $i)$ we cannot define a denumerable sequence $f'(x_{1}, 0),$ $f'(x_{1}, 1), \ldots$ such that:

\vspace{+1ex}
\begin{quote}
(e) $f(k, i) = f'(k, i)$ for all $i \geq 0$.
\end{quote}

\vspace{+1ex}
We cannot thus define a deterministic Turing machine $TM$ that will give evidence that the PA formula $[(\exists_{1} x_{3})F(x_{1}, x_{2}, x_{3})]$ interprets as true under the standard interpretation of PA for any given sequence of numerals $[(a_{1}, a_{2})]$.

\vspace{+.5ex}
Hence $[(\exists_{1} x_{3})F(x_{1}, x_{2}, x_{3})]$ is not algorithmically computable over $\mathbb{N}$ under the standard interpretation of PA.
\end{quote}

\noindent The theorem follows.\hfill $\Box$

\vspace{+1ex}
The above theorem now suggests the following definition:

\begin{definition}
\label{sec:2.3.def.1}
\noindent \textbf{Effective computability:} A number-theoretic function is effectively computable if, and only if, it is algorithmically verifiable.
\end{definition}

\vspace{+1ex}
However, we still need to formally define what it means for a number-theoretic formula of an arithmetical language to be:

\begin{quote}
(i) Algorithmically verifiable;

\vspace{+.5ex}
(ii) Algorithmically computable.
\end{quote}

\noindent under an interpretation.

\section{Interpretation of an arithmetical language in terms of Turing computability}
\label{sec:1.03}

We begin by noting that we can, in principle, define the classical `satisfaction' and `truth' of the formulas of a first order arithmetical language, such as PA, \textit{verifiably} under an interpretation using as \textit{evidence} the computations of a simple functional language.

Such definitions follow straightforwardly for the atomic formulas of the language (i.e., those without the logical constants that correspond to `negation', `conjunction', `implication' and `quantification') from the standard definition of a deterministic Turing machine\footnote{\cite{Me64}, pp.229-231.}, based essentially on Alan Turing's seminal 1936 paper on computable numbers\footnote{\cite{Tu36}.}.

Moreover, it follows from Alfred Tarski's seminal 1933 paper on the the concept of truth in the languages of the deductive sciences\footnote{\cite{Ta33}.} that the `satisfaction' and `truth' of those formulas of a first-order language which contain logical constants can be inductively defined, under an interpretation, in terms of the `satisfaction' and `truth' of the interpretations of only the atomic formulas of the language.

Hence the `satisfaction' and `truth' of those formulas (of an arithmetical language) which contain logical constants can, in principle, also be defined verifiably under an interpretation using as evidence the computations of a deterministic Turing machine.

We show in Section \ref{sec:5.4.0} that this is indeed the case for PA under the standard interpretation $\mathcal{I}_{PA(\mathbb{N},\ Standard)}$, when this is explicitly defined as in Section \ref{sec:5.4.0.0}.

We show in Section \ref{sec:6.2}, moreover, that we can further define `algorithmic truth' and `algorithmic falsehood' under $\mathcal{I}_{PA(\mathbb{N},\ Standard)}$ such that the PA axioms interpret as always algorithmically true.

\begin{quote}
\footnotesize
\textbf{Significance of `algorithmic truth':} The \textit{algorithmically} true propositions of $\mathbb{N}$ under $\mathcal{I}_{PA(\mathbb{N},\ Standard)}$ are thus a proper subset of the \textit{verifiably} true propositions of $\mathbb{N}$ under $\mathcal{I}_{PA(\mathbb{N},\ Standard)}$; and suggest a possible finitary `model'\footnote{\cite{Me64}, p.51.} of PA that would establish the consistency of PA constructively.
\end{quote}

\subsection{The definitions of `algorithmic truth' and `algorithmic falsehood' under $\mathcal{I}_{PA(\mathbb{N},\ Standard)}$ are not symmetric with respect to `truth' and `falsehood' under $\mathcal{I}_{PA(\mathbb{N},\ Standard)}$}
\label{sec:1.03.1}

However, the definitions (in Section \ref{sec:6.2}) of `algorithmic truth' and `algorithmic falsehood' under $\mathcal{I}_{PA(\mathbb{N},\ Standard)}$ are not symmetric with respect to classical (verifiable) `truth' and `falsehood' under $\mathcal{I}_{PA(\mathbb{N},\ Standard)}$.

For instance, if a formula $[F(x_{1}, x_{2}, \ldots, x_{n})]$ of an arithmetic is algorithmically true under an interpretation (such as $\mathcal{I}_{PA(\mathbb{N},\ Standard)}$) that appeals to the evidence provided by the computations of a deterministic Turing machine, then, for any given denumerable sequence of numerical values $[a_{1}, a_{2}, \ldots]$, the formula $[F(a_{1}, a_{2}, \ldots, a_{n})]$ is also algorithmically true under the interpretation.

\begin{quote}
\footnotesize
\textbf{Denumerable sequence:} We shall presume that any such sequence is `given' in the sense of being defined by the `evidence' of a deterministic Turing machine.
\end{quote}

In other words, there is a deterministic Turing machine which can provide evidence that the interpretation $F^{*}(a_{1}, a_{2}, \ldots, a_{n})$ of $[F(a_{1}, a_{2}, \ldots, a_{n})]$ \textit{holds} in $\mathbb{N}$ for any given denumerable sequence of natural numbers $(a_{1}, a_{2}, \ldots)$.

\begin{quote}
\footnotesize
\textbf{Defining the term `hold':} We define the term `hold'---when used in connection with an interpretation of a formal language and, more specifically, with reference to the operations of a deterministic Turing machine associated with the atomic formulas of the language---explicitly in Section \ref{sec:5.4.0}; the aim being to avoid appealing to the classically subjective (and existential) connotation implicitly associated with the term under an implicitly defined standard interpretation of an arithmetic\footnote{As, for instance, in \cite{Go31}.}.
\end{quote}

However, if a formula $[F(x_{1}, x_{2}, \ldots, x_{n})]$ of an arithmetic is algorithmically false under an interpretation that appeals to the evidence provided by the computations of a deterministic Turing machine, we cannot conclude that there is a denumerable sequence of numerical values $[a_{1}, a_{2}, \ldots]$ such that the formula $[F(a_{1}, a_{2}, \ldots, a_{n})]$ is algorithmically false under the interpretation.

Reason: If a formula $[F(x_{1}, x_{2}, \ldots, x_{n})]$ of an arithmetic is algorithmically false under such an interpretation, then we can only conclude that there is no deterministic Turing machine which can provide evidence that the interpretation $F^{*}(x_{1}, x_{2}, \ldots,$ $x_{n})$ holds in $\mathbb{N}$ for any given denumerable sequence of natural numbers $(a_{1}, a_{2}, \ldots)$; we cannot conclude that there is a denumerable sequence of natural numbers $(b_{1}, b_{2}, \ldots)$ such that $F^{*}(b_{1}, b_{2}, \ldots, b_{n})$ does not hold in $\mathbb{N}$.

Such a conclusion would require:

\vspace{+1ex}
(i) either some additional evidence that will verify for some assignment of numerical values to the free variables of $[F]$ that the corresponding interpretation $F^{*}$ does not hold\footnote{Essentially reflecting Brouwer's objection to the assumption of Aristotle's particularisation over an infinite domain.};

\vspace{+1ex}
(ii) or the additional assumption that either Aristotle's particularisation holds over the domain of the interpretation (as is implicitly presumed under the standard interpretation of PA) or the arithmetic is $\omega$-consistent\footnote{An assumption explicitly introduced by G\"{o}del in \cite{Go31}.}.

\begin{quote}
\footnotesize
\textbf{An issue of consistency:} An issue that we do not address here is whether the assumption of Aristotle's particularisation (or that of $\omega$-completeness) is consistent with the evidence provided by the computations of a deterministic Turing machine for defining the satisfaction and truth of the formulas of an arithmetic under an interpretation.   
\end{quote}

\section{Defining algorithmic verifiability and algorithmic computability}
\label{sec:1.02}

The asymmetry of Section \ref{sec:1.03.1} suggests the following two concepts\footnote{My thanks to Dr. Chaitanya H. Mehta for advising that the focus of this investigation should be the distinction between these two concepts.}: 

\begin{definition}
\label{sec:1.02.def.1}
\textbf{Algorithmic verifiability}: An arithmetical formula $[F(x_{1}, x_{2},$ $\ldots, x_{n})]$ is algorithmically verifiable under an interpretation if, and only if, for any given sequence of numerals $[a_{1}, a_{2}, \ldots, a_{n}]$, we can define a deterministic Turing machine $TM$ that computes $[F(a_{1}, a_{2}, \ldots,$ $a_{n})]$ and will halt on null input if, and only if, $[F(a_{1}, a_{2}, \ldots, a_{n})]$ interprets as either true or false under the interpretation.
\end{definition}

\begin{quote}
\scriptsize
\textbf{Interpretation of an arithmetical language \textit{verifiably} in terms of Turing computability} Of course such a definition requires defining the `algorithmic verifiability' of the formulas of an arithmetical language under an interpretation in terms of the operations of a deterministic Turing machine. Such definitions follow straightforwardly for the atomic formulas of the language (i.e., those without the logical constants that correspond to `negation', `conjunction', `implication' and `quantification') from Alan Turing's seminal 1936 paper on computable numbers\footnote{\cite{Tu36}.}. Moreover, it follows from Alfred Tarski's seminal 1933 paper on the the concept of truth in the languages of the deductive sciences\footnote{\cite{Ta33}.} that the `algorithmic verifiability' of the formulas of a formal language which contain logical constants can be inductively defined under an interpretation in terms of the `algorithmic verifiability' of the interpretations of the atomic formulas of the language. Hence the `algorithmic verifiability' of the formulas containing logical constants in the above definition can, in principle, be defined in terms of the operations of a deterministic Turing machine. We show in Section \ref{sec:5.4.0} that this is indeed the case.
\end{quote}

\begin{definition}
\label{sec:1.02.def.2}
\textbf{Algorithmic computability}: An arithmetical formula $[F(x_{1},$ $x_{2}, \ldots,$ $x_{n})]$ is algorithmically computable under an interpretation if, and only if, we can define a deterministic Turing machine $TM$ that computes $[F(x_{1},$ $x_{2}, \ldots,$ $x_{n})]$ and, for any given sequence of numerals $[a_{1}, a_{2}, \ldots, a_{n}]$, will halt on $[a_{1}, a_{2}, \ldots, a_{n}]$ if, and only if, $[F(a_{1}, a_{2}, \ldots, a_{n})]$ interprets as either true or false under the interpretation.
\end{definition}

\begin{quote}
\scriptsize
\textbf{Interpretation of an arithmetical language \textit{algorithmically} in terms of Turing computability} We show in Section \ref{sec:5.4.0} that the `algorithmic computability' of the formulas of a formal language which contain logical constants can also be inductively defined under an interpretation in terms of the `algorithmic computability' of the interpretations of the atomic formulas of the language.
\end{quote}

We now show that the above concepts are well-defined under the standard interpretation of PA.

\section{Standard definitions of the `satisfiability' and `truth' of a formal language under an interpretation}
\label{sec:5.4.0}

We note that standard interpretations of the formal reasoning and conclusions of classical first order theory---based primarily on the work of Cantor, G\"{o}del, Tarski, and Turing---seem to admit the impression that the classical, Tarskian, truth (satisfiability) of the propositions of a formal mathematical language under an interpretation is, both, non-algorithmic and essentially unverifiable constructively\footnote{Thus giving rise to J.\ R.\ Lucas' G\"{o}delian Argument in \cite{Lu61}; but see also \cite{An07a}, \cite{An07b}, \cite{An07c} and \cite{An08}.}.

However---if mathematics is to serve as a universal set of languages of, both, precise expression and unambiguous communication---such interpretations may need to be balanced by an alternative, constructive and intuitionistically unobjectionable, interpretation---of classical foundational concepts---in which non-algorithmic truth (satisfiability) is defined effectively\footnote{Such an interpretation is broadly outlined in \cite{An07d}}.

For instance, we note that---essentially following standard expositions\footnote{cf. \cite{Me64}, p.51.} of Tarski's inductive definitions on the `satisfiability' and `truth' of the formulas of a formal language under an interpretation---we can define:

\begin{definition}
\label{sec:5.def.1}
If $[A]$ is an atomic formula $[A(x_{1}, x_{2}, \ldots, x_{n})]$ of a formal language S, then the denumerable sequence $(a_{1}, a_{2}, \ldots)$ in the domain $\mathbb{D}$ of an interpretation $\mathcal{I}_{S(\mathbb{D})}$ of S satisfies $[A]$ if, and only if:

\begin{quote}
(i) $[A(x_{1}, x_{2}, \ldots, x_{n})]$ interprets under $\mathcal{I}_{S(\mathbb{D})}$ as a relation $A^{*}(x_{1}, x_{2},$ $\ldots, x_{n})$ in $\mathbb{D}$ for a witness $\mathcal{W}_{\mathbb{D}}$ of $\mathbb{D}$;

\vspace{+1ex}
(ii) there is a Satisfaction Method, SM($\mathcal{I}_{S(\mathbb{D})}$) that provides objective \textit{evidence}\footnote{In the sense of \cite{Mu91}.} by which any witness $\mathcal{W}_{\mathbb{D}}$ of $\mathbb{D}$ can \textbf{define} for any atomic formula $[A(x_{1}, x_{2}, \ldots, x_{n})]$ of S, and any given denumerable sequence $(b_{1}, b_{2}, \ldots)$ of $\mathbb{D}$, whether the proposition $A^{*}(b_{1}, b_{2}, \ldots, b_{n})$ holds or not in $\mathbb{D}$;

\vspace{+1ex}
(iii) $A^{*}(a_{1}, a_{2}, \ldots, a_{n})$ holds in $\mathbb{D}$ for $\mathcal{W}_{\mathbb{D}}$.
\end{quote}
\end{definition}

\begin {quote}
\footnotesize
\textbf{Witness:} From a constructive perspective, the existence of a `witness' as in (i) above is implicit in the usual expositions of Tarski's definitions.

\vspace{+1ex}
\textbf{Satisfaction Method:} From a constructive perspective, the existence of a Satisfaction Method as in (ii) above is also implicit in the usual expositions of Tarski's definitions.

\vspace{+1ex}
\textbf{A constructive perspective:} We highlight the word `\textit{\textbf{define}}' in (ii) above to emphasise the constructive perspective underlying this paper; which is that the concepts of `satisfaction' and `truth' under an interpretation are to be explicitly viewed as objective assignments by a convention that is witness-independent. A Platonist perspective would substitute `decide' for `define', thus implicitly suggesting that these concepts can `exist', in the sense of needing to be discovered by some witness-dependent means---eerily akin to a `revelation'---even when the domain $\mathbb{D}$ is $\mathbb{N}$.
\end{quote}

We can now inductively assign truth values of `satisfaction', `truth', and `falsity' to the compound formulas of  a first-order theory S under the interpretation $\mathcal{I}_{S(\mathbb{D})}$ in terms of \textit{only} the satisfiability of the atomic formulas of S over $\mathbb{D}$ as follows\footnote{Compare \cite{Me64}, p.51; \cite{Mu91}.}:

\begin{definition}
\label{sec:5.def.2}
A denumerable sequence $s$ of $\mathbb{D}$ satisfies $[\neg A]$ under $\mathcal{I}_{S(\mathbb{D})}$ if, and only if, $s$ does not satisfy $[A]$;
\end{definition}

\begin{definition}
\label{sec:5.def.3}
A denumerable sequence $s$ of $\mathbb{D}$ satisfies $[A \rightarrow B]$ under $\mathcal{I}_{S(\mathbb{D})}$ if, and only if, either it is not the case that $s$ satisfies $[A]$, or $s$ satisfies $[B]$;
\end{definition}

\begin{definition}
\label{sec:5.def.4}
A denumerable sequence $s$ of $\mathbb{D}$ satisfies $[(\forall x_{i})A]$ under $\mathcal{I}_{S(\mathbb{D})}$ if, and only if, given any denumerable sequence $t$ of $\mathbb{D}$ which differs from $s$ in at most the $i$'th component, $t$ satisfies $[A]$;
\end{definition}

\begin{definition}
\label{sec:5.def.5}
A well-formed formula $[A]$ of $\mathbb{D}$ is true under $\mathcal{I}_{S(\mathbb{D})}$ if, and only if, given any denumerable sequence $t$ of $\mathbb{D}$, $t$ satisfies $[A]$;
\end{definition}

\begin{definition}
\label{sec:5.def.6}
A well-formed formula $[A]$ of $\mathbb{D}$ is false under under $\mathcal{I}_{S(\mathbb{D})}$ if, and only if, it is not the case that $\mathbb{D}$ is true under $\mathcal{I}_{S(\mathbb{D})}$.
\end{definition}

It follows that\footnote{cf. \cite{Me64}, pp.51-53.}:

\begin{theorem}
\label{sec:5.thm.1}
(\textit{Satisfaction Theorem}) If, for any interpretation $\mathcal{I}_{S(\mathbb{D})}$ of a first-order theory S, there is a Satisfaction Method SM($\mathcal{I}_{S(\mathbb{D})}$) which holds for a witness $\mathcal{W}_{\mathbb{D}}$ of $\mathbb{D}$, then:

\begin{quote}
(i) The $\Delta_{0}$ formulas of S are decidable as either true or false over $\mathbb{D}$ under $\mathcal{I}_{S(\mathbb{D})}$;

\vspace{+1ex}
(ii) If the $\Delta_{n}$ formulas of S are decidable as either true or as false over $\mathbb{D}$ under $\mathcal{I}_{S(\mathbb{D})}$, then so are the $\Delta(n+1)$ formulas of S.
\end{quote}
\end{theorem}

\vspace{+1ex}
\noindent \textbf{Proof} It follows from the above definitions that:

\vspace{+1ex}
(a) If, for any given atomic formula $[A(x_{1}, x_{2}, \ldots, x_{n})]$ of S, it is decidable by $\mathcal{W}_{\mathbb{D}}$ whether or not a given denumerable sequence $(a_{1}, a_{2}, \ldots)$ of $\mathbb{D}$ satisfies $[A(x_{1}, x_{2}, \ldots, x_{n})]$ in $\mathbb{D}$ under $\mathcal{I}_{S(\mathbb{D})}$ then, for any given compound formula $[A^{1}(x_{1}, x_{2}, \ldots, x_{n})]$ of S containing any one of the logical constants $\neg, \rightarrow, \forall$, it is decidable by $\mathcal{W}_{\mathbb{D}}$ whether or not $(a_{1}, a_{2}, \ldots)$ satisfies $[A^{1}(x_{1}, x_{2}, \ldots, x_{n})]$ in $\mathbb{D}$ under $\mathcal{I}_{S(\mathbb{D})}$;

\vspace{+1ex}
(b) If, for any given compound formula $[B^{n}(x_{1}, x_{2}, \ldots, x_{n})]$ of S containing $n$ of the logical constants $\neg, \rightarrow, \forall$, it is decidable by $\mathcal{W}_{\mathbb{D}}$ whether or not a given denumerable sequence $(a_{1}, a_{2}, \ldots)$ of $\mathbb{D}$ satisfies $[B^{n}(x_{1}, x_{2}, \ldots, x_{n})]$ in $\mathbb{D}$ under $\mathcal{I}_{S(\mathbb{D})}$ then, for any given compound formula $[B^{(n+1)}(x_{1}, x_{2}, \ldots, x_{n})]$ of S containing $n+1$ of the logical constants $\neg, \rightarrow, \forall$, it is decidable by $\mathcal{W}_{\mathbb{D}}$ whether or not $(a_{1}, a_{2}, \ldots)$ satisfies $[B^{(n+1)}(x_{1}, x_{2}, \ldots, x_{n})]$ in $\mathbb{D}$ under $\mathcal{I}_{S(\mathbb{D})}$;

\vspace{+1ex}
We thus have that:

\vspace{+1ex}
(c) The $\Delta_{0}$ formulas of S are decidable by $\mathcal{W}_{\mathbb{D}}$ as either true or false over $\mathbb{D}$ under $\mathcal{I}_{S(\mathbb{D})}$;

\vspace{+1ex}
(d) If the $\Delta_{n}$ formulas of S are decidable by $\mathcal{W}_{\mathbb{D}}$ as either true or as false over $\mathbb{D}$ under $\mathcal{I}_{S(\mathbb{D})}$, then so are the $\Delta(n+1)$ formulas of S. \hfill $\Box$

\vspace{+1ex}
In other words, if the atomic formulas of of S interpret under $\mathcal{I}_{S(\mathbb{D})}$ as decidable with respect to the Satisfaction Method SM($\mathcal{I}_{S(\mathbb{D})}$) by a witness $\mathcal{W}_{\mathbb{D}}$ over some domain $\mathbb{D}$, then the propositions of S (i.e., the $\Pi_{n}$ and $\Sigma_{n}$ formulas of S) also interpret as decidable with respect to SM($\mathcal{I}_{S(\mathbb{D})}$) by the witness $\mathcal{W}_{\mathbb{D}}$ over $\mathbb{D}$.

\subsection{An explicit definition of the standard interpretation of PA}
\label{sec:5.4.0.0}

We now consider the application of Tarski's definitions to the standard interpretation $\mathcal{I}_{PA(\mathbb{N},\ Standard)}$ of PA where:

\begin{quote}
(a) we define S as PA with standard first-order predicate calculus as the underlying logic;

\vspace{+1ex}
(b) we define $\mathbb{D}$ as $\mathbb{N}$;

\vspace{+1ex}
(c) we take SM($\mathcal{I}_{PA(\mathbb{N},\ Standard)}$) as a deterministic Turing machine.
\end{quote}

We note that:

\vspace{+1ex}
\begin{theorem}
\label{sec:5.5.lem.1}
The atomic formulas of PA are algorithmically verifiable under the standard interpretation $\mathcal{I}_{PA(\mathbb{N},\ Standard)}$. 
\end{theorem}

\vspace{+1ex}
\textbf{Proof} If $[A(x_{1}, x_{2}, \ldots, x_{n})]$ is an atomic formula of PA then, for any given denumerable sequence of numerals $[b_{1}, b_{2}, \ldots]$, the PA formula $[A(b_{1}, b_{2},$ $\ldots, b_{n})]$ is an atomic formula of the form $[c=d]$, where $[c]$ and $[d]$ are atomic PA formulas that denote PA numerals. Since $[c]$ and $[d]$ are recursively defined formulas in the language of PA, it follows from a standard result\footnote{For any natural numbers $m,\ n$, if $m \neq n$, then PA proves $[\neg(m = n)]$ (\cite{Me64}, p.110, Proposition 3.6). The converse is obviously true.} that, if PA is consistent, then $[c=d]$ interprets as the proposition $c=d$ which either holds or not for a witness $\mathcal{W}_{\mathbb{N}}$ in $\mathbb{N}$.

Hence, if PA is consistent, then $[A(x_{1}, x_{2}, \ldots, x_{n})]$ is algorithmically verifiable since, for any given denumerable sequence of numerals $[b_{1}, b_{2}, \ldots]$, we can define a deterministic Turing machine $TM$ that will compute $[A(b_{1}, b_{2}, \ldots, b_{n})]$ and halt on null input as evidence that the PA formula $[A(b_{1}, b_{2}, \ldots, b_{n})]$ is decidable under the interpretation.

\vspace{+1ex}
The theorem follows.\hfill $\Box$

\vspace{+1ex}
We thus have that:

\begin{corollary}
The `satisfaction' and `truth' of PA formulas containing logical constants can be defined under the standard interpretation of PA in terms of the evidence provided by the computations of a deterministic Turing machine.
\end{corollary}

\subsection{Defining `algorithmic truth' under the standard interpretation of PA and interpreting the PA axioms}
\label{sec:6.2}

Now we note that, in addition to Theorem \ref{sec:5.5.lem.1}:

\vspace{+1ex}
\begin{theorem}
\label{sec:5.5.lem.2}
The atomic formulas of PA are algorithmically computable under the standard interpretation $\mathcal{I}_{PA(\mathbb{N},\ Standard)}$. 
\end{theorem}

\vspace{+1ex}
\textbf{Proof} If $[A(x_{1}, x_{2}, \ldots, x_{n})]$ is an atomic formula of PA then we can define a deterministic Turing machine $TM$ that will compute $[A(x_{1}, x_{2}, \ldots, x_{n})]$ and halt on any given denumerable sequence of numerals $[b_{1}, b_{2}, \ldots]$ as evidence that the PA formula $[A(b_{1}, b_{2}, \ldots, b_{n})]$ is decidable under the interpretation.

\vspace{+1ex}
The theorem follows.\hfill $\Box$

This suggests the following definitions:

\begin{definition}
\label{sec:6.2.def.1}
A well-formed formula $[A]$ of PA is algorithmically true under $\mathcal{I}_{PA(\mathbb{N},\ Standard)}$ if, and only if, there is a deterministic Turing machine which computes $[A]$ and provides evidence that, given any denumerable sequence $t$ of $\mathbb{N}$, $t$ satisfies $[A]$;
\end{definition}

\begin{definition}
\label{sec:6.2.def.2}
A well-formed formula $[A]$ of PA is algorithmically false under under $\mathcal{I}_{PA(\mathbb{N})}$ if, and only if, it is not algorithmically true under $\mathcal{I}_{PA(\mathbb{N})}$.
\end{definition}

The significance of defining `algorithmic truth' under $\mathcal{I}_{PA(\mathbb{N},\ Standard)}$ as above is that:

\begin{lemma}
\label{sec:6.2.lem.1}
The PA axioms PA$_{1}$ to PA$_{8}$ are algorithmically computable as algorithmically true over $\mathbb{N}$ under the interpretation $\mathcal{I}_{PA(\mathbb{N},\ Standard)}$.
\end{lemma}

\vspace{+1ex}
\noindent \textbf{Proof} Since $[x+y]$, $[x \star y]$, $[x = y]$, $[{x^{\prime}}]$ are defined recursively\footnote{cf. \cite{Go31}, p.17.}, the PA axioms PA$_{1}$ to PA$_{8}$ interpret as recursive relations that do not involve any quantification. The lemma follows straightforwardly from Definitions \ref{sec:5.def.1} to \ref{sec:5.def.6} in Section \ref{sec:5.4.0} and Theorem \ref{sec:5.5.lem.1}.\hfill $\Box$

\begin{lemma}
\label{sec:6.2.lem.2}
For any given PA formula $[F(x)]$, the Induction axiom schema $[F(0)$ $\rightarrow (((\forall x)(F(x) \rightarrow F(x^{\prime}))) \rightarrow (\forall x)F(x))]$ interprets as algorithmically true under $\mathcal{I}_{PA(\mathbb{N},\ Standard)}$.
\end{lemma}

\vspace{+1ex}
\noindent \textbf{Proof} By Definitions \ref{sec:5.def.1} to \ref{sec:6.2.def.2}:

\begin{quote}

(a) If $[F(0)]$ interprets as algorithmically false under $\mathcal{I}_{PA(\mathbb{N},\ Standard)}$ the lemma is proved.

\vspace{+1ex}
\begin{quote}
\footnotesize Since $[F(0) \rightarrow (((\forall x)(F(x) \rightarrow F(x^{\prime}))) \rightarrow (\forall x)F(x))]$ interprets as algorithmically true if, and only if, either $[F(0)]$ interprets as algorithmically false or $[((\forall x)(F(x) \rightarrow F(x^{\prime}))) \rightarrow (\forall x)F(x)]$ interprets as algorithmically true.
\end{quote}

\vspace{+1ex}
(b) If $[F(0)]$ interprets as algorithmically true and $[(\forall x)(F(x) \rightarrow F(x^{\prime}))]$ interprets as algorithmically false under $\mathcal{I}_{PA(\mathbb{N},\ Standard)}$, the lemma is proved.

\vspace{+1ex}
(c) If $[F(0)]$ and $[(\forall x)(F(x) \rightarrow F(x^{\prime}))]$ both interpret as algorithmically true under $\mathcal{I}_{PA(\mathbb{N},\ Standard)}$, then by Definition \ref{sec:6.2.def.1} there is a deterministic Turing machine $TM$ that, for any natural number $n$, will give evidence that the formula $[F(n) \rightarrow F(n^{\prime})]$ is true under $\mathcal{I}_{PA(\mathbb{N},\ Standard)}$.

\vspace{+1ex}
Since $[F(0)]$ interprets as algorithmically true under $\mathcal{I}_{PA(\mathbb{N},\ Standard)}$, it follows that there is a deterministic Turing machine $TM$ that, for any natural number $n$, will give evidence that the formula $[F(n)]$ is true under the interpretation.

\vspace{+1ex}
Hence $[(\forall x)F(x)]$ is algorithmically true under $\mathcal{I}_{PA(\mathbb{N},\ Standard)}$.
\end{quote}

Since the above cases are exhaustive, the lemma follows.\hfill $\Box$

\begin{quote}
\footnotesize
\textbf{The Poincar\'{e}-Hilbert debate:} We note that Lemma \ref{sec:6.2.lem.2} appears to settle the Poincar\'{e}-Hilbert debate\footnote{See \cite{Hi27}, p.472; also \cite{Br13}, p.59; \cite{We27}, p.482; \cite{Pa71}, p.502-503.} in the latter's favour. Poincar\'{e} believed that the Induction Axiom could not be justified finitarily, as any such argument would necessarily need to appeal to infinite induction. Hilbert believed that a finitary proof of the consistency of PA was possible. 
\end{quote}

\begin{lemma}
\label{sec:6.2.lem.3}
Generalisation preserves algorithmic truth under $\mathcal{I}_{PA(\mathbb{N},\ Standard)}$.
\end{lemma}

\vspace{+1ex}
\noindent \textbf{Proof} The two meta-assertions:

\begin{quote}
`$[F(x)]$ interprets as algorithmically true under $\mathcal{I}_{PA(\mathbb{N},\ Standard)}$\footnote{See Definition \ref{sec:5.def.5}}'

\vspace{+.5ex}
and

\vspace{+.5ex}
`$[(\forall x)F(x)]$ interprets as algorithmically true under $\mathcal{I}_{PA(\mathbb{N},\ Standard)}$'
\end{quote}

\noindent both mean:

\begin{quote}
$[F(x)]$ is algorithmically computable as always true under $\mathcal{I}_{PA(\mathbb{N}),}$ $_{Standard)}$. \hfill $\Box$
\end{quote} 

\vspace{+1ex}
It is also straightforward to see that:

\begin{lemma}
\label{sec:6.2.lem.4}
Modus Ponens preserves algorithmic truth under $\mathcal{I}_{PA(\mathbb{N},\ Standard)}$. \hfill $\Box$
\end{lemma}

We thus have that:

\begin{theorem}
\label{sec:6.2.lem.5}
The axioms of PA are always algorithmically true under the interpretation $\mathcal{I}_{PA(\mathbb{N},\ Standard)}$, and the rules of inference of PA preserve the properties of algorithmic satisfaction/truth under $\mathcal{I}_{PA(\mathbb{N},\ Standard)}$, without appeal to Aristotle's particularisation.\hfill $\Box$
\end{theorem}

\section{Appendix A: G\"{o}del's Theorem VII(2)}
\label{sec:5.4.1}
\begin{center}
\footnotesize{(Excerpted from \cite{Go31} pp.29-31.)}
\end{center}

\vspace{+1ex}
Every relation of the form $x_{0} = \phi(x_{1}, \ldots, x_{n})$, where $\phi$ is recursive, is arithmetical and we apply complete induction on the rank of $\phi$. Let $\phi$ have rank $s (s>1)$. \ldots

\vspace{+1ex}
$\phi(0, x_{2}, \ldots, x_{n}) = \psi(x_{2}, \ldots, x_{n})$

$\phi(k+1, x_{2}, \ldots, x_{n}) = \mu[k, \phi(k, x_{2}, \ldots, x_{n}), x_{2}, \ldots, x_{n}]$

\noindent(where $\psi, \mu$ have lower rank than $s$).

\vspace{+1ex}
\ldots we apply the following procedure: one can express the relation $x_{0} = \phi(x_{1}, \ldots, x_{n})$ with the help of the concept ``sequence of numbers" $(f)$\footnote{$f$ denotes here a variable whose domain is the sequence of natural numbers. The $(k+1)$st term of a sequence $f$ is designated $f_{k}$ (and the first, $f_{0}$).} in the following manner:

\begin{quote}
$x_{0} = \phi(x_{1}, \ldots, x_{n}) \sim (\exists f)\{f_{0} = \psi(x_{2}, \ldots, x_{n})\ \&\ (\forall k)(k<x_{1} \rightarrow f_{k+1} = \mu(k, f_{k}, x_{2}, \ldots, x_{n})\ \&\ x_{0} = f_{x_{1}}\}$
\end{quote}

If $S(y, x_{2}, \ldots, x_{n}), T(z, x_{1}, \ldots, x_{n+1})$ are the arithmetical relations\\ which, according to the inductive hypothesis, are equivalent to $y = \psi(x_{2},$ $\ldots, x_{n})$, and $z = \mu(x_{1}, \ldots, x_{n+1})$ respectively, then we have:

\begin{quote}
$x_{0} = \phi(x_{1}, \ldots, x_{n}) \sim (\exists f)\{S(f_{0}, x_{2}, \ldots, x_{n})\ \&\ (\forall k)[k<x_{1} \rightarrow$ \\ $T(f_{k+1}, k, x_{2}, \ldots, x_{n})]\ \&\ x_{0} = f_{x_{1}}\}$ \hspace{+29ex} (17)
\end{quote}

Now we replace the concept ``sequence of numbers" by ``pairs of numbers" by correlating with the number pair $n, d$ the sequence of numbers $f^{(n, d)}$ ($f_{k}^{(n, d)} = [n]_{1+(k+1)d}$, where $[n]_{p}$ denotes the smallest non-negative remainder of $n$ modulo $p$).

Then:

\vspace{+1ex}
Lemma 1: If $f$ is an arbitrary sequence of natural numbers and $k$ is an arbitrary natural number, then there exists a pair of natural numbers $n, d$ such that $f^{(n, d)}$ and $f$ coincide in their first $k$ terms.

\vspace{+1ex}
Proof: Let $l$ be the greatest of the numbers $k, f_{0}, f_{1}, \ldots, f_{k-1}$. Determine $n$ so that

\vspace{+1ex}
$n \equiv f_{i}\ [mod\ (1+(i+1)l!)]$ for $i = 0, 1, \ldots, k-1$,

\vspace{+1ex}
\noindent which is possible, since any two of the numbers $1+(i+1)l!\ (i = 0, 1, \ldots, k-1)$ are relatively prime. For, a prime dividing two of these numbers must also divide the difference $(i_{1} - i_{2})l!$ and therefore, since $i_{1} - i_{2} < l$, must also divide $l!$, which is impossible. The number pair $n, l!$ fulfills our requirement.

\vspace{+1ex}
Since the relation $x = [n]_{p}$ is defined by

\vspace{+1ex}
$x \equiv n\ (mod\ p)\ \&\ x < p$

\vspace{+1ex}
\noindent and is therefore arithmetical, then so also is the relation $P(x_{0}, x_{1}, \ldots, x_{n})$ defined as follows:

\begin{quote}
$P(x_{0}, x_{1}, \dots, x_{n}) \equiv (\exists n, d)\{S([n]_{d+1}, x_{2}, \ldots, x_{n})\ \&\ (\forall k)[k < x_{1} \rightarrow T([n]_{1+d}$ $_{(k+2)}, k,$ $[n]_{1+d(k+1)}, x_{2}, \ldots, x_{n})]\ \&\ x_{0} = [n]_{1+d(x_{1} + 1)}\}$
\end{quote}

\noindent which, according to (17) and Lemma 1, is equivalent to $x_{0} = \phi(x_{1}, \ldots, x_{n})$ (in the sequence $f$ in (17) only its values up to the $(x+1)$th term matter). Thus, Theorem VII(2) is proved.

\begin{quote}
\footnotesize{
\textbf{Comment}: G\"{o}del's remark that ``in the sequence $f$ in (17) only its values up to the $(x+1)$th term matter" is significant. The proof of Theorem \ref{sec:2.3.thm.1} depends upon the fact that the equivalence between $f^{(n, d)}$ and $f$ cannot be extended non-terminatingly.
} 
\end{quote}

\noindent \tiny{Authors postal address: 32 Agarwal House, D Road, Churchgate, Mumbai - 400 020, Maharashtra, India.\ Email: re@alixcomsi.com, anandb@vsnl.com.}

\end{document}